\documentclass[11pt]{amsart}

\usepackage[margin=1in]{geometry}
\usepackage{setspace}
\usepackage{amsmath,amssymb}
\usepackage{booktabs}
\usepackage{tikz}
\usepackage{pgfplots}
\usepackage{float}
\usepackage{hyperref}
\usepackage{makecell}

\usepackage{tikz}
\usepackage{pgfplots}
\pgfplotsset{compat=1.16}

\usepackage{amsthm}

\title{Reducing the Incentive to Tank:  The Ex Post Gold Plan}

\author{Bret Benesh}
\address{Department of Mathematics, College of Saint Benedict and Saint John's University\\37 College Avenue South\\ Saint Joseph, MN 56374\\ bbenesh@csbsju.edu}

\date{\today}

\begin{document}

\maketitle

\begin{abstract}

Many recent proposals for reducing tanking in draft lotteries share a common structure: losses improve draft position early in the season while wins improve draft position later. While such systems improve late-season incentives, they retain a predictable pivot point that tanking teams can exploit strategically. This paper proposes a simple modification that introduces uncertainty into the timing of the incentive switch. The proposed metric, the \emph{Realized Elimination Wins Determinant} (REWIND), ranks teams according to the number of wins obtained after their ex post elimination date, which makes this a variation of the Gold Plan. Because the ex post elimination date cannot be known with certainty during the season, the mechanism weakens incentives for strategic losing while preserving incentives for competitive effort after elimination.  Moreover, the ex post elimination date is typically earlier than other proposed pivot points, so there is a longer period where a tanking team's best strategy is to win.  The Ex Post Gold plan uses the REWIND metric to create a simple system where every team will be incentivized to win at least half of their games in most seasons.

\end{abstract}

\section{Introduction}

Draft lotteries attempt to balance competitive by helping losing teams become competitive. However, systems that reward poor performance risk encouraging strategic losing, commonly known as \emph{tanking}.  A growing class of reform proposals attempts to address this issue by dividing the season into two phases with a pivot at time $T_p$:

\begin{itemize}
\item losses improve draft position before time $T_p$.
\item wins improve draft position after time $T_p$.
\end{itemize}

Examples of potential $T_p$ include the All-Star break, trade deadline, March $1$st, and the $60$-game mark.  The current policy defines $T_p$ to be the end of the season, so wins never improve a team's lottery position.  

While such mechanisms improve late-season incentives, they share an important limitation: the pivot time $T_p$ is known to teams in advance. As a result, a tanking team can confidently lose before $T_P$ and try to win after $T_P$ to maximize their lottery standing.  This paper proposes a simple modification that introduces uncertainty into the timing of the incentive change $T_P$, making it impossible for a team to know whether losing or winning a particular game will improve their lottery standing for a portion of the season.

We will be using the National Basketball Association (NBA) as a context throughout this paper.  The sixth seed is the last team in each conference that is guaranteed to be in the playoffs, so we will use the sixth seed as a reference when we discuss being eliminated from the playoffs.

\section{The Gold Plan and Literature Review}
Lenten \cite{lenten2016} argues that anti-tanking reform should resolve draft-order incentives earlier in the season, and Lenten, Smith, and Boys \cite{lenten2018} evaluate an AFL policy that orders teams by elimination from finals contention rather than by end-of-season standings. These papers show the appeal of elimination-based reform, but they also leave open a timing problem: when the pivot is mathematically or administratively observable, teams can often infer when losses stop helping. Banchio and Munro \cite{banchiomunro2021} prove that any draft system based only on final performance with unequal weights creates incentives for teams to underperform, and then constructs a way of creating a dynamic stopping time for each team based on observed performance.  Their mechanism is mathematically sophisticated and dynamically adjusts draft incentives based on observed performance, but may be cahllenging to implement or communicate in practice.

The Gold Plan \cite{gold2010}, implemented in the Professional Women’s Hockey League, is a variation of plans stated above that uses each team's date of mathematical elimiation.  That is, there is not a single $T_p$ that applies to every team, but rather each team $i$ has its own time $T_i$ that acts as a pivot from losses counting toward draft position to wins counting toward draft position.   The Gold Plan also deviates from the plans above by only counting wins after $T_i$, but not losses prior to $T_i$.  The compensation for not counting losses is that worse teams will experience  mathematical elimination earlier, which gives them more chances to accumulate wins after $T_i$.   For the purposes

The Gold Plan improves incentives after elimination, but the tanking teams are still incentivized to lose prior to mathematical elimiation and win after.

\section{The Ex Post Gold Plan}
A slight modification of the Gold Plan dramatically increases a tanking team's incentive to win.    We define the Ex Post Gold Plan below using the date when a team was ex post eliminated from the playoffs.  

First, we define the \emph{Realized Elimination Wins Determinant} (REWIND) for a team $i$ to be the number of wins after a date $T_i$, where $T_i$ is determined by the following algorithm.

\begin{enumerate}
\item At the end of the season, determine the number of losses $L$ of the sixth seed from $i$'s conference.
\item Define $T_i$ to be the date when $i$ loses their $L+1$st game.  Call $T_i$ the \emph{REWIND Date} and $L+1$ the \emph{REWIND Target}.
\item Define $W_i$ to be the total number of wins for $i$ at the end of the season, and define $R_i$, the \emph{REWIND wins}, to be the total number of wins for $i$ at $T_i$.
\item The \emph{REWIND Score for $i$} is $W_i-R_i$. 
\end{enumerate}

Note that $T_i$, the ex post elimination date, is exactly the date when team $i$ would not have gotten the sixth seed (or higher), even if they had won every game after $T_i$.   Also, note the difficulty of predicting $T_i$ prior to the end of the season:  a tanking team needs to predict how many losses the sixth seed will have, and they do not even know who the sixth seed will be at the end of the season.  The ex post nature of $T_i$ is what makes it difficult for the tanking team to strategize.

Note that the ex post elimination date is consistent with what we do with the mathematical elimination in the usual Gold Plan.  Even the standard notion of mathematical elimination depends on how one models feasible future outcomes. In particular, simple win-based calculations may overstate a team's chances by ignoring interactions among competing teams.  For example, suppose that teams $A$, $B$, and $C$ are competing for one playoff spot, all teams have two games left, and $A$ has two fewer wins than both $B$ and $C$.  Team $A$ is not technically mathematically eliminated in the usual sense, since they could tie $B$ and $C$ if $A$ wins its remaining two games and $B$ and $C$ both lose their remaining two games.  However, this will be impossible if $B$ and $C$ play each other in one of the two remaining games.  Thus, the usual notion of ``mathematical elimination" is overly simplified at times, and the REWIND Score mirrors this simplicity.

This algorithm reduces the incentive to tank for part of the season by tying the pivot date $T_i$ to future outcomes.  This makes it impossible to know $T_i$ with certainty when it happens, but $T_i$ is easy to compute at the end of the season prior to determining the lottery order.

\section{Optimal Tanking Strategy}

We now consider how a tanking team would strategize to maximize the chances of lottery success. The Ex Post Gold Plan uses the number of losses for the team in the sixth seed, which has had a minimum of 28 losses (the 2007--2008 Utah Jazz \cite{bref2008}) and a maximum of 41 losses (most recently from the 2014--2015 Milwaukee Bucks \cite{bref2015}) in the 82-game seasons since 2000.  Then a team that wants to tank under the Ex Post Gold Plan would likely consider the season to be in three parts.

\begin{enumerate}
\item A tanking team should definitely lose in the first part of the season, which occurs historically between the times when the team has their $0$th loss and their $29$th loss.
\item A tanking team not be certain about whether to win or lose in the second part of the season, which occurs historically between their $30$th loss and $41$nd loss. 
\item A tanking team should definitely try to win in the third part of the season, which occurs historically after the $42$nd loss.
\end{enumerate}

There is a lot of uncertainty in the second part, which spans $12$ games---about $15\%$ of the season.  This will be illustrated in an example below.  This would require good and careful modeling about how the rest of the season will play out \emph{for the entire league}, which will be difficult to do.  

Moreover, the third section of the season of the season---where the tanking team is definitely incentivized to win---is historically equivalent to almost half of the season since each REWIND Target this century has been at most 42.  That is, it would be a rare season when any team would want to win fewer than 50\% of their games.

Consider an example from the 2023--2024 season, where the Portland Trailblazers had a conference-worst 21--67 record.  Figure~\ref{fig:POR2024} compares how their lottery standing would have been determined under the current policy, the Gold Plan, and the Ex Post Gold Plan.    Here is how the Ex Post Gold Plan was determined in detail.

\begin{enumerate}
\item At the end of the season, the sixth place team was the New Orleans Pelicans, who ended up with 33 losses for the entire season. 
\item We define $T_{POR2024}$ to be the ex post elimination date February 2, 2024, since that is the day they lost their 34th game of the season.
\item We define $W_{POR2024}$ to be 15, since Portland was 15--34 on February 2, 2024.  
\item Portland ended the season with 21 wins, so their RECORD Score is $21-15=6$. 
\end{enumerate}

\begin{table}[H]
\centering
\begin{tabular}{ccccccc}
\toprule
\makecell{Anti-Tanking \\ Strategy} & Pivot Date & \makecell{POR Record \\ at Pivot Date} &  \makecell{Relevant Team and\\ Record} & \makecell{Reason for \\ Pivot} & \makecell{Metric} \\
\midrule
Current NBA Plan & End of Season & 21--61& NA& Definition &  $61$ \\
Gold Plan & 3/11/2024 & 18--46& \makecell{PHO at Pivot \\(38--27)}& $46+38 > 82$ &  $21-18=3$ \\
\makecell{\\Ex Post Gold Plan} & 2/2/2024 & 15--34& \makecell{NOP at end of \\ season (49--33)}& $34 > 33$ & $21-15=6$ \\
\bottomrule
\end{tabular}
\caption{\label{fig:POR2024}Comparison of Calculations Lottery Plans for the Portland Trailblazers 2023--2024.  Their record was 21--61.}
\end{table}


To explain the table further:  in the Gold Plan, the reason why Portland was not mathematically eliminated on March 10th is that one could conceive of a universe where Portland wins all of its last 19 games and Phoenix loses all of its 18 games.The Ex Post Gold Plan does not allow for such a universe, though, since the benefit of hindsight allows us to see that Phoenix actually finished with 49 wins.  In fact, Phoenix didn't even get the 6th seed---the New Orleans Pelicans did (due to a tie-breaker with Phoenix). 

Under the Ex Post Gold Plan, the only universe that we care about is the one that was actually realized---the one where the New Orleans Pelicans finished in sixth place with a 49--33 record, which is why the REWIND Target is 34 (this was also the REWIND Target for every other Western Conference team that year).  With hindsight, we can see that Portland was actually eliminated earlier than March 11th:  they got their 34th loss on February 2nd, which means that---again, with hindsight---they had no chance of clinching a guaranteed playoff position on February 2nd, their REWIND Date.  

Now we will explain some of the complexities a tanking team needs to consider in the second part of the season.  There was no way for Portland to know that their REWIND Date was going to be February 2nd.  As a small example, New Orleans lost to the Los Angeles Lakers on the last game of the season.  Had they won, they would have finished 50--32, and Portland's REWIND Date would have been January 28th instead of February 2nd.  

Consider Portland's game on January 29th.  As it turns out, the optimal tanking strategy for the Ex Post Gold Plan was to lose in order to get to the REWIND Date faster.  However, if the New Orleans Pelicans had beaten the Lakers on April 14th, it would have turned out that the optimal tanking strategy would have been to win that game since the hypothetical REWIND Date of January 28th had passed.  Determining the optimal tanking strategy for this particular game requires knowledge of the future, which is a strength of the Ex Post Gold Plan.

See the Appendix for an evaluation of how this would affect the lottery teams for the post-Covid 82-game NBA seasons (2021--2022 to 2024--2025). For instance, the Utah Jazz and the Philadelphia 76ers were often cited as teams that tanked in 2024--2025, and Figure~\ref{fig:SingleUniverseGold2025} shows that they would have dropped six and nine places in the lottery, respectively.

\section{Varitations, Advantages, and Disadvantages of the Ex Post Gold Plan}

The Ex Post Gold Plan can be altered by changing the REWIND Target to be based off the number of losses of any of the 7th through 10th seeds.  In order to introduce more uncertainty, a lottery could be held at the end of the season to determine which of the 6th, 7th, 8th, 9th, or 10th seeds will be used to determine the REWIND Target.  

The advantages of the Ex Post Gold Plan include the following.

\begin{enumerate}
\item  Teams are incentivized to win for more of the season.
\item Like the usual Gold Plan, losses never directly help the lottery rankings.
\item The optimal tanking  strategy for any one game requires knowledge of the future.  
\item This model is simple:  at the end of the season, you just count the number of wins after they have gotten to their REWIND Target.
\item The tanking team that intentionally loses a midseason game risks \emph{ex post regret}; if the realized elimination date ultimately places the game in the portion of the season where wins improve draft position, the loss represents a foregone opportunity to improve the team's lottery standing.  This could be particularly strong with coaches who have strong influence over the front office---the coach wants to win, and it would be particularly painful to intentionally lose a game when winning turned out to be the better tanking strategy.
\end{enumerate}

There are a couple of known disadvantages to the Ex Post Gold Plan.

\begin{enumerate}
\item It is harder for truly bad teams to rise in the lottery, since wins are the only thing that determines lottery standing.  The current NBA lottery rules could be kept to guarantee that the team with the most losses picks no lower than fifth---REWIND would be used to determine the lottery odds, the lottery would be used to determine the teams that draft in the first four picks, and the total number of losses could be used to determine the rest of the lottery.  
\item A truly great team could get the first pick if they have injuries in the first half of the season.  This is a problem under the current rules, too (David Robinson's injury led to the Spurs drafting Tim Duncan), but it would be exacerbated under the Ex Post Gold Plan.  Tying the Ex Post Gold Plan to, say, the 10th seed instead of the sixth seed would help counteract this fear.  
\end{enumerate}

\section{Conclusion}

The Ex Post Gold Plan incentive is a variation of plans that use wins to improve lottery standing after a certain date.  This variation is an improvement via two mechanisms.

\begin{enumerate}
\item It extends the time when a tanking team is incentivized to win by using the ex post elimination date, which occurs earlier in the season that most of the other proposals.  It does this in a non-arbitrary way, since the ex post elimination date has intrinsic meaning.
\item It reduces the incentive to win for the period immediately before the ex post elimination date, since the ex post elimination date is unknown until the end of the season.
\end{enumerate}

The Ex Post Gold Plan is simple to understand and makes it harder for a tanking team to identify when they should lose games.  When tied to the sixth seed, this will normally provide incentive for every team in the league to try to win at least half of their games.

\section{Appendix}

Below are examples of how the lottery would have gone in the last four seasons under the Ex Post Gold Plan.  All tie-breakers were determined by total losses.

\begin{table}[H]
\centering
\begin{tabular}{cccccccc|c}
\toprule
\makecell{Actual \\Rank} & Team & Record &  \makecell{REWIND \\Date} & \makecell{REWIND\\ Wins} & \makecell{REWIND \\Score} & \makecell{REWIND\\ Rank} &\makecell{Rank\\ Change} & \makecell{REWIND\\ Ranking}\\
\midrule
1 & HOU & 20--62& 1/28& 14 & 6 & 7 &-6 &  ORL\\
2 & ORL & 22--60& 1/12& 7& 15 & 1 & +1& SAS\\
3 & DET & 23--59& 1/23& 11& 12 &3&+0& DET\\
4 & OKC & 24--58& 2/5& 17& 7 &6&-2&NYK\\
5 & IND & 25--57& 2/4& 19& 6 &8&-3&SAC \\
6 & POR  & 27--55& 2/24 & 25 & 2 &13&-7& OKC\\
7 & SAC & 30--52& 2/3 & 19& 11 &5&+2& HOU\\
8 & LAL & 33--49& 3/3 & 27& 6 &9&-1&IND\\
9 & SAS &34--48 &2/9 &20 & 14 &2&+7&LAL\\
10 & WAS & 35--47& 3/9& 29& 6 &10&+0&WAS\\
11 & NYK & 37--45& 3/2& 25& 12 &4&+7&LAC\\
12 & LAC & 42--40&3/14 & 36& 6 &11&+1&CHO\\
13 & CHO & 43--39 & 3/27 & 39 & 4 &12&+1&POR\\
14 & CLE & 44--38 & 4/5&43 & 1 &14&+0&CLE\\
\bottomrule
\end{tabular}
\caption{REWIND Results 2021--2022 \cite{bref2022e,bref2022w}.  The sixth seeds were the Chicago Bulls at 46--36 and the Denver Nuggets at 48--34, so the REWIND Targets were 37 for the East and 35 for the West.}
\end{table}

\begin{table}[H]
\centering
\begin{tabular}{cccccccc|c}
\toprule
\makecell{Actual \\Rank} & Team & Record &  \makecell{REWIND \\Date} & \makecell{REWIND\\ Wins} & \makecell{REWIND \\Score} & \makecell{REWIND\\ Rank} &\makecell{Rank\\ Change} & \makecell{REWIND\\ Ranking}\\
\midrule
1 & DET & 17--35& 1/28& 13 & 4 & 7 &-6 &  CHO\\
2 & SAS & 22--60& 2/3& 14& 8 & 3 & -1& HOU\\
3 & HOU & 22--60& 2/3& 13& 9 &2&+1& SAS\\
4 & CHO & 27--55& 2/2& 15& 12 &1&+3&ORL\\
5 & POR & 33--49& 3/17& 31& 2 &11&-6&CHI\\
6 & ORL  & 34--48& 3/5 & 27 & 7 &4&+2& TOR\\
7 & WAS & 35--47& 3/17 & 32& 3 &9&-2& DET\\
8 & IND & 35--47& 3/13 & 31& 4 &8&+0&IND\\
9 & UTA &37--45 &3/25 &35 & 2 &12&-3&WAS\\
10 & DAL & 38--44& 3/26& 36& 2 &13&-3&OKC\\
11 & OKC & 40--42& 3/28& 37& 3 &10&+1&POR\\
12 & CHI & 40--42&3/22 & 34& 6 &5&+7&UTA\\
13 & TOR & 41--41 & 3/22 & 35 & 6 &6&+7&DAL\\
14 & NOP & 42--40 & 4/4&40 & 2 &14&+0&NOP\\
\bottomrule
\end{tabular}
\caption{REWIND Results 2022--2023 \cite{bref2023e,bref2023w}.  The sixth seeds were the Brooklyn Nets at 45--37 and the Los Angeles Clippers at 44--38, so the REWIND Targets were 38 for the East and 39 for the West.}
\end{table}

\begin{table}[H]
\centering
\begin{tabular}{cccccccc|c}
\toprule
\makecell{Actual \\Rank} & Team & Record &  \makecell{REWIND \\Date} & \makecell{REWIND\\ Wins} & \makecell{REWIND \\Score} & \makecell{REWIND\\ Rank} &\makecell{Rank\\ Change} & \makecell{REWIND\\ Ranking}\\
\midrule
1 & DET & 14--68& 1/12& 3 & 11 & 4 &-3 & TOR\\
2 & WAS & 15--67& 1/24& 7& 8 & 9 & -7& HOU\\
3 & CHO & 21--61& 1/31& 10& 11 &5&-2& SAS\\
4 & POR & 21--61& 2/2& 15& 6 &11&-7&DET\\
5 & SAS & 22--60& 1/19& 7& 15 &3&+2&CHO\\
6 & TOR  & 25--57& 2/14 & 19 & 6 &1&+5& BRK\\
7 & MEM & 27--55& 2/8 & 18& 9 &8&-1& GSW\\
8 & UTA & 31--51& 3/2 & 27& 4 &13&-5&MEM\\
9 & BRK  &32--50 &2/27 &22 & 10 &6&+3&WAS\\
10 & ATL & 36--46& 3/13& 29& 7 &10&+0&ATL\\
11 & CHI & 39--43& 3/21& 34& 5 &12&-1&POR\\
12 & HOU & 41--41&2/29 & 25& 16 &2&+10&CHI\\
13 & SAC & 46--36 & 4/9 & 45 & 1 &14&-1&UTA\\
14 & GSW & 46--36 & 3/23&36 & 10 &7&+7&SAC\\
\bottomrule
\end{tabular}
\caption{REWIND Results 2023--2024 \cite{bref2024e,bref2024w}.  The sixth seeds were the Orlando Magic at 47--35 and the New Orleans Pelicans at 49--33, so the REWIND Targets were 36 for the East and 34 for the West.}
\end{table}

\begin{table}[H]
\centering
\begin{tabular}{cccccccc|c}
\toprule
\makecell{Actual \\Rank} & Team & Record &  \makecell{REWIND \\Date} & \makecell{REWIND\\ Wins} & \makecell{REWIND \\Score} & \makecell{REWIND\\ Rank} &\makecell{Rank\\ Change} & \makecell{REWIND\\ Ranking}\\
\midrule
1 & UTA & 17--65& 1/27& 10 & 7 & 7 &-6 & TOR\\
2 & WAS & 18--64& 1/26& 6& 12 & 2 & +0& WAS\\
3 & CHO & 19--63& 2/12& 13& 6 &9&-6& CHI \\
4 & NOP & 21--61& 1/25& 12& 9 &4&+0& NOP\\
5 & PHI  & 24--58& 3/3& 21& 3 &14&-9&SAS\\
6 & BRK  & 26--56& 3/1 & 21 & 5 &11&-5& POR\\
7 & TOR & 30--52& 2/21 & 17& 13 &1&+6& UTA\\
8 & SAS & 34--48& 3/2 & 25& 9 &5&+3&PHO\\
9 & PHO  &36--46 &3/7 &29 & 7 &8&+1&CHO\\
10 & POR & 36--46& 3/2& 27& 9 &6&+4&DAL\\
11 & DAL & 39--43& 3/12& 33& 6 &10&+1&BRK\\
12 & CHI & 39--43&3/15 & 28& 11 &3&+9&SAC\\
13 & ATL & 40--42 & 4/1 & 36 & 4 &13&+0&ATL\\
14 & SAC & 40--42 & 3/20&35 & 5 &12&+2&PHI\\
\bottomrule
\end{tabular}
\caption{\label{fig:SingleUniverseGold2025}REWIND Results 2024--2025 \cite{bref2025e,bref2025w}.  The sixth seeds were the Detroit Pistons at 44--38 and the Minnesota Timberwolves at 49--33, so the REWIND Targets were 39 for the East and 34 for the West.}
\end{table}

\end{document}